\documentclass[12pt]{article}
\usepackage{amssymb}
\usepackage{latexsym}

\newtheorem{thm}{Theorem}[section]
\newtheorem{prop}[thm]{Proposition} 
 \newtheorem{dfn}[thm]{Definition}

\newcommand{\reals}{{\mathbb R}}

\newcommand{\ad}{{\rm ad}}

\newcommand{\cala}{{\cal A}}

\newcommand{\cale}{{\cal E}}
\newcommand{\calf}{{\cal F}}

\newcommand{\call}{{\cal L}}

\newcommand{\half}{\textstyle{\frac{1}{2}}}
\newcommand{\third}{\textstyle{\frac{1}{3}}}

\newcommand{\arrows}{\,\lower1pt\hbox{$\longrightarrow$}\hskip-.24in\raise2pt
             \hbox{$\longrightarrow$}\,}
\newcommand{\bleft}{[\![}
\newcommand{\bright}{]\!]}
\newcommand{\smalcirc}{\mbox{\tiny{$\circ $}}}

\newcommand{\glnr}{{\mathfrak gl}(n,\reals)}
\title{{\bf Omni-Lie Algebras
}
}
\author{Alan Weinstein\thanks{Research partially supported by NSF
Grant DMS-99-71505}
\\Department of Mathematics\\
University of California\\
Berkeley, CA 94720 USA\\
{\small(alanw@math.berkeley.edu)}}
\date{}
\begin{document}
\maketitle
\section{Introduction}
\label{sec:intro}

We introduce on the space $\cale_n ={\mathfrak gl}(n,\reals)\times \reals^n$
the antisymmetric bracket operation 
\begin{equation}
  \label{eq:bracket}
  \bleft (A_1,v_1),(A_2,v_2)\bright =
  ([A_1,A_2],\textstyle{\frac{1}{2}}(A_1 v_2 -A_2 v_1)).
\end{equation}
Without the factor of $\frac{1}{2}$, this would be the
semidirect product Lie algebra for the usual action of $\glnr$ on
$\reals^n$.  With the factor of $\frac{1}{2}$, the bracket does not satisfy the
Jacobi identity.  Nevertheless, it does satisfy the Jacobi identity on
many subspaces which are closed under the bracket.  In fact, we will
see that any Lie algebra structure on $\reals^n$ is realized on such a
subspace.  

If $B$ is any bilinear operation on $\reals^n$, we define the
adjoint operator $\ad_B:\reals^n\rightarrow \glnr$ by
$\ad_B(v)(w)=B(v,w)$, and we denote by $\calf_B \subset \cale_n $ the graph
of this operator.  A simple calculation shows:
\begin{prop}
  \label{prop:realize}
If $B$ is a skew symmetric bilinear operation on $\reals^n$, then $B$
satisfies the Jacobi identity if and only if $\calf_B$ is closed under the
bracket operation (\ref{eq:bracket}) on $\cale_n $.  When this condition is satisfied, the
restriction to $\calf_B$ of the natural projection from $\cale_n  $ to
$\reals^n$ is an isomorphism between the restricted $\cale_n  $ bracket
and the operation $B$.  
\end{prop}

The skew symmetry of an operation $B$ can also be seen as a
property of the graph $\calf_B$.  We introduce a symmetric bilinear
form on $\cale_n $ with values in $\reals^n$:
\begin{equation}
  \label{eq:symmetric}
  \langle (A_1,v_1),(A_2,v_2)\rangle=\half(A_1 v_2 + A_2 v_1).
\end{equation}
The operation $B$ is skew symmetric if and only if $\calf_B$ is
isotropic for this bilinear form.  When this is the case, $\calf_B$ is
actually a maximal isotropic subspace.

As we mentioned above, the bracket operation $\bleft ~,~\bright$ does
not satisfy the Jacobi identity.  In fact, we have a simple formula
for the jacobiator 
\begin{equation}
  \label{eq:jacobiator}
  J(e_1,e_2,e_3)=\bleft \bleft e_1,e_2\bright,e_3\bright + \mbox{c.p.},
\end{equation}
where the $e_j$ are elements of $\cale_n $ and ``c.p.'' means to add the
two terms obtained by cyclic permutation of the three indices.
In terms of the ``Cartan 3-form'' 
\begin{equation}
  \label{eq:cartan}
  T(e_1,e_2,e_3)=\textstyle{\frac{1}{3}}\langle\bleft e_1,e_2\bright,
  e_3\rangle +\mbox{c.p.},
\end{equation}
the jacobiator is given by:
\begin{equation}
  \label{eq:jacobianomaly}
  J(e_1,e_2,e_3)=(0,T(e_1,e_2,e_3)).
\end{equation}

Proposition \ref{prop:realize} follows immediately from Equation
\ref{eq:jacobianomaly} and the fact that the graph $\calf_B$ is {\em
maximal} isotropic.  

We define a $D$-structure on $\reals^n$ to be any maximal isotropic
subspace of $\cale_n$ which is closed under the bracket operation.  By
(\ref{eq:jacobianomaly}), any $D$-structure is a Lie algebra with the
restricted bracket.  Among the $D$-structures are not only the
$n$-dimensional Lie algebras, realized on the graphs of their adjoint
representations, but also the ``horizontal'' subspace $\glnr \oplus
\{0\}$.  Unfortunately, this is the only graph of a mapping from
$\glnr$ to $\reals^n$ which is isotropic in $\cale_n$.  The
reader is invited to
find other $D$-structures.

We come now to the principal question raised by the construction
above.  All $n$-dimensional Lie algebras can be embedded in the
space $\cale_n$, which has a bracket operation which does not quite
satisfy the Jacobi identity.  We will refer to $\cale_n$ as an {\bf
omni-Lie algebra}.  {\em Is there a global object
corresponding to $\cale_n$, obtained by some kind of ``integration'',
which would serve as an ``omni-Lie group''?}

\noindent
{\bf Remark}.  In response to the posting of a preprint version of
this work, Michael Kinyon made some suggestions which have led to the
resolution of some (but not all!) of the problems raised here.
Details will appear in \cite{ki-we:double}.

\section{Courant algebroids}
\label{sec:courant}

Although our presentation of the omni-Lie algebras in the previous
section was self contained, in fact we came to this construction by
linearizing at a point the following bracket on
the sections of $E=TM\oplus T^*M$, where $M$ is a differentiable
manifold.  It was introduced by T.~Courant \cite{co:dirac}.
\begin{equation}
  \label{eq:courant}
  \bleft (\xi_1,\theta_1),(\xi_2,\theta_2)\bright =
  ([\xi_1,\xi_2],\call_{\xi_1}\theta_2 -\call_{\xi_2}\theta_1 -\half
  d(i_{\xi_1}\theta_2  -i_{\xi_2}\theta_1)),
\end{equation}
where $\call_{\xi}$ and $i_{\xi}$ are the operations of Lie derivative
and interior product by the vector field $\xi$.  

Courant introduced his bracket to unify the treatment of Poisson
structures (bivector fields $\pi$ on $M$ for which the bracket
$\{f,g\}=\pi(df,dg)$ on $C^{\infty}(M)$ satisfies the Jacobi identity)
and presymplectic structures (2-forms on $M$ which are closed). 
Each bivector field $\pi$ or 2-form $\omega$ on $M$ gives rise to a graph which
is a subbundle $F$ of $E$ whose fibres are maximal isotropic with respect
the the (indefinite) inner product
\begin{equation}
  \label{eq:inner}
  \langle(\xi_1,\theta_1),(\xi_2,\theta_2)\rangle =
  \half(\theta_2(\xi_1)+\theta_1(\xi_2)). 
\end{equation}
The space of sections of the subbundle $F$ is closed under the {\bf
  Courant bracket}
(\ref{eq:courant}) if and only if $\pi$ [resp. $\omega$] is a Poisson
[resp. presymplectic] structure.  Courant defined a
{\bf Dirac structure} on $M$ to be any maximal isotropic subbundle
$F\subset TM\oplus T^*M$ whose sections are closed under the Courant
bracket.  Dirac structures include not only Poisson and
presymplectic structures, but also foliations on $M$.  If $B$ is an
integrable subbundle of $TM$, $B^{\perp}$ its annihilator in $T^*M$,
then $B\oplus B^{\perp}$ is a Dirac structure.  

Any Dirac structure $F$ on $M$ is a {\bf Lie algebroid}.  This means that:
(1) the sections of $M$ have a Lie algebra structure $[~,~]$ (over $\reals$);
(2) there is a bundle map $\rho:F\rightarrow TM$, called the {\bf anchor},
which induces a Lie algebra homomorphism from sections of $F$ to
vector fields on $M$; (3) for any sections $e_1$ and $e_2$ of $F$ and
a function $f$ on $M$, 
\begin{equation}
\label{eq:anchor}
[e_1,fe_2]=f[e_1,e_2]+(\rho(e_1)\cdot f)e_2.
\end{equation}
For a Dirac structure, the bracket and anchor are the
restriction to $F$ of the Courant bracket and the projection on the
first summand of $TM\oplus T^*M$.  When the Dirac structure is a
Poisson structure $\pi$, the subbundle $F$ may be identified with $T^*M$ by
projection onto the second summand of $E$, and the resulting Lie
algebroid structure on $T^*M$ is the infinitesimal object
 corresponding to (local) symplectic groupoids
 \cite{co-da-we:groupoides}\cite{va:lectures} for the Poisson manifold $(M,\pi)$.

Thus, the bundle $E$ carries a structure which does not quite satisfy
the axioms of a Lie algebroid, since its bracket does not satisfy the
Jacobi identity, but it contains many subbundles on which the
restricted bracket {\em is} a Lie algebroid structure; these include
the Lie algebroids of the symplectic groupoids of all the Poisson
structures on $M$.  Since Lie
algebroids are the infinitesimal objects for (local) Lie groupoids, it
is natural to ask whether there is a global, groupoid-like object
corresponding to $E$ which contains all these groupoids.

The properties of Courant's bracket were the basis for the definition
of a {\bf Courant algebroid} by Liu, Xu, and the author
\cite{li-we-xu:manin}.  This object is defined to be a vector bundle
$E$ over a manifold $M$ carrying a field of inner products (i.e.
nondegenerate symmetric
bilinear forms) along the fibres, an antisymmetric bracket operation
$\bleft~,~\bright$ on its space of sections, and a bundle map $\rho:E\rightarrow
TM$ such that the Jacobi identity and the Leibniz rule
(\ref{eq:anchor}) are satisfied modulo terms which are differentials,
in a certain sense, of terms involving the inner products.  In
addition, the bilinear form itself satisfies some conditions, one of
which is a modified version of ``adjoint invariance''.  (The precise
axioms can be found in algebraic form in Section \ref{sec:algebra}
below.)   When $M$ is a
point, all the error terms vanish, and a Courant algebroid is just a Lie
algebra with an adjoint-invariant inner product.  In this case, the
corresponding global object is clearly a Lie group with a bi-invariant
(possibly indefinite) metric.  

In any Courant algebroid, one may define the Dirac structures to be
the maximal isotropic subbundles whose sections are closed under the
bracket; since the anomalies vanish on Dirac structures, they are Lie
algebroids.  When $M$ is a point, a pair of complementary
Dirac structures is a
Manin triple corresponding to a Lie bialgebra, and conversely the
direct sum of a Lie bialgebra and its dual is a Courant algebroid over
a point.  This is the double of the Lie bialgebra, and the global
object is the double of a Poisson Lie group corresponding to the Lie
bialgebra.  
Similarly, there is a notion of Lie bialgebroid due to
Mackenzie and Xu \cite{ma-xu:lie}.  The double of a Lie bialgebroid is
a Courant algebroid \cite{li-we-xu:manin}; what is missing is the
double of the Poisson groupoid corresponding to the Lie bialgebroid.

Some progress has been made in relating Courant algebroids to other
algebraic structures.  In \cite{ro-we:courant}, it is shown that
Courant algebroids can be considered as strongly homotopy Lie algebras.
In Roytenberg's thesis \cite{ro:courant}, an approach to Lie algebroids in terms of
homological vector fields on supermanifolds was developed to describe
arbitrary Courant algebroids.  The thesis also develops the idea,
suggested by some calculations in \cite{li-we-xu:manin} and
observations by Y.~Kosmann-Schwarzbach and P.~Severa, that a
non-antisymmetric version of the bracket on a Courant algebra is an
example of the Leibniz algebras introduced by Loday \cite{lo:version}
and called Loday algebras in \cite{ko:from}, where
Kosmann-Schwarzbach shows how Loday brackets can be obtained as
so-called derived brackets from Poisson brackets and derivations.  It
is this version of the Courant bracket which plays the central role in
Roytenberg's work.  

Although P.~Severa has observed that a class of Courant algebroids
obtained by deforming Courant's original example may be seen as the
infinitesimal objects corresponding to gerbes (see \cite{br:loop} for
a discussion of gerbes), there is still no satisfactory description of
the groupoid-like object corresponding to a general Courant algebroid,
nor of the group-like object corresponding to the sections of a
Courant algebroid.  In the hope of clarifying the situation, we may
try to ``linearize'' Courant's original example at a point of
$M=\reals^n$.  The result, as we shall explain next, is the omni-Lie
algebra $\cale_n$ of Section \ref{sec:intro}.

\section{$C$-Algebras}
\label{sec:algebra}
There is an algebraic version of the notion of Lie algebroid, in which
the role of a vector bundle over a manifold is played by the
algebraic analogue of its space of sections, namely a module over a
commutative algebra.  This concept goes under various names, including
``Lie-Rinehart algebra'' and $(R,\cala)$ Lie algebra; we refer
the reader to \cite{ma:survey} for an extensive list of them.
In this section, we introduce an 
analogous algebraic version of Courant algebroids.

We begin with a commutative ground ring $R$ and a commutative
$R$-algebra $\cala$ (not necessarily unital).  Next, we consider an
$\cala$-module $\cale$ with a homomorphism $\rho$ from $\cale$ to the
$\cala$-module of derivations of $\cala$.  For any $f\in A$, we define
its $\cale$ differential $d_{\cale }f$ in the dual $\cala$ module
$\cale^*$ by $d_{\cale}f (e)=\rho(e)f.$ Now suppose that $\cale$ also
carries a symmetric $\cala$-bilinear form $\langle ~,~ \rangle$ which
is {\bf weakly nondegenerate} in the sense that the associated
homomorphism $\beta$ from $\cale$ to the dual $\cala$-module $\cale^*$
is injective.  If $d_{\cale}f$ is in the image of $\beta$ for a given 
$f\in \cala$, we define the {\bf gradient} $Df$ to be the well-defined element 
$\half\beta^{-1}d_{\cale}f$ of $\cale$; i.e. 
$\langle D f , e\rangle = \half  \rho (e) f $ for all $e\in \cale$.
(The factor of $\half$ is
irritating, but it has to go somewhere; we have chosen to put it
here to match the conventions in \cite{li-we-xu:manin}.)

\begin{dfn}
  \label{dfn:calgebra}
With the definitions and notation above, an $(R,\cala)$ $C$-{\bf algebra}
 is an $\cala$-module $\cale$
carrying a nondegenerate $\cala$-valued symmetric bilinear form
$\langle~,~\rangle$, an antisymmetric $R$-bilinear operation
$\bleft~,~\bright$, and an $\cala$-module homomorphism $\rho$ with values
in the $R$-derivations of $\cala$ such that the following properties are
satisfied:
\begin{enumerate}
\setcounter{enumi}{-1}
\item The gradient $Df$ is defined for all $f\in \cala$;
\item For any $e_{1}, e_{2}, e_{3}$ in $\cale$, 
$\bleft\bleft e_{1}, e_{2}\bright , e_{3}\bright +\mbox{c.p.}=
D  T(e_{1}, e_{2}, e_{3});$
\item  for any $e_{1}, e_{2} $ in $\cale$,
$\rho \bleft e_{1}, e_{2}\bright =[\rho e_{1}, \rho  e_{2}];$
\item  for any $e_{1}, e_{2} $ in $ \cale$ and $f$ in $ \cala$,
$\bleft e_{1}, fe_{2}\bright =f\bleft e_{1}, e_{2}\bright +(\rho
(e_{1})f)e_{2}- \langle e_{1}, e_{2}\rangle D f ;$
\item $\rho \smalcirc D =0$, i.e.,  for any $f, g$ in $\cala$
$\langle D f,  D  g\rangle=0$;
\item for any $e, h_{1}, h_{2} $ in $\cale$,
  $\rho (e) \langle h_{1}, h_{2}\rangle =\langle \bleft e, h_{1}\bright +D \langle e
  ,h_{1}\rangle  , 
h_{2}\rangle +\langle h_{1}, \bleft e , h_{2}\bright +D  \langle e
,h_{2}\rangle \rangle $,
\end{enumerate}
where $T(e_{1}, e_{2}, e_{3})$ is the element of $\cala$
defined by:
\begin{equation}
\label{eq:T0}
 T(e_{1}, e_{2}, e_{3})=\third \langle \bleft e_{1}, e_{2} \bright , e_{3}\rangle +\mbox{c.p.},
\end{equation}
\end{dfn}

When $R=\reals$ and $\cala$ is the algebra of smooth functions on a
manifold $M$, an $(R,\cala)$ $C$-algebra is just the space of sections
of a Courant algebroid.  On the other hand, to to see the omni-Lie
algebras of Section \ref{sec:intro} as $(R,\cala)$ $C$-algebras, we
let $R$ be $\reals$ and $\cala$ be $\reals^n$ with the multiplication
in which all products are zero.  Geometrically, we should think of the
latter as the algebra of functions which are defined on an
infinitesimal neighborhood of the origin in the dual space of
$\reals^n$ and which vanish at the origin.  (One could adjoin the
constant functions at the cost of complicating the example slightly.)
$\cale_n$ is $\glnr\oplus \reals^n$ with the module structure in which
all scalar products are zero.  We identify the derivations of $\cala$
with $\glnr$, so that $\rho$ can be defined as projection on the first
summand of $\glnr\oplus \reals^n.$   The inner product and bracket are
given by (\ref{eq:symmetric}) and (\ref{eq:bracket}).  The gradient D
is then given for all $v$ by $Dv=(0,v)$.  The axioms of an
$(R,\cala)$ $C$-algebra may be checked directly, or as a
consequence of identities satisfied in the original Courant algebroid.

It is clear that isotropic subalgebras of $(R,\cala)$ $C$-algebras are Lie
algebras.  The fact that {\em all} $n$-dimensional Lie algebras arise
this way in the omni-Lie algebra $\glnr\oplus \reals ^n$ can be seen
as a consequence of the fact that, among the Dirac structures on the
dual space of $\reals^n$ are the Lie-Poisson structures attached to
all Lie algebra structures on $\reals^n$.  

Finally, we note that Evens and Lu
\cite{ev-lu:variety} have recently shown that the variety of maximal
isotropic subalgebras in the double of a Lie bialgebra
carries a natural Poisson
structure.  It would be interesting to see whether their work
extends to our more general setting.


\begin{thebibliography}{99}
\bibitem{br:loop}
Brylinski, J.L., {\em Loop spaces, characteristic classes and geometric
quantization}. Birkh\"auser, Boston, 1993.
\bibitem{co-da-we:groupoides}
Coste, A., Dazord, P., et Weinstein, A., Groupo\"{\i}des symplectiques,
{\em Publications du D\'{e}partement de Math\'{e}matiques,
Universit\'{e} Claude Bernard--Lyon I} {\bf 2A} (1987), 1--62.  
\bibitem{co:dirac} 
Courant, T.J., Dirac manifolds, {\em Trans. A.M.S.} {\bf 319} (1990),
631--661. 
\bibitem{ev-lu:variety}
Evens, S., and Lu, J-H., On the variety of Lagrangian subalgebras,
preprint math.DG/9909005. 
\bibitem{ki-we:double}
Kinyon, M.K., and Weinstein, A.,
A double construction for Courant algebroids (in preparation).
\bibitem{ko:from}
Kosmann-Schwarzbach, Y., From Poisson algebras to Gerstenhaber algebras,
{\em Ann. Inst. Fourier (Grenoble)} {\bf 46} (1996), 1243--1274.
\bibitem{li-we-xu:manin}
Liu, Z.-J., Weinstein, A., and Xu, P., Manin triples for Lie
bialgebroids, {\em J. Diff. Geom} {\bf 45} (1997), 547--574.
\bibitem{lo:version}
Loday, J.L., Une version non commutative des alg\`ebres de Lie: les
alg\`ebres de 
Leibniz, {\em Enseign. Math.} {\bf 39} (1993), 269--293. 
\bibitem{ma:survey}
Mackenzie, K.C.H.,
Lie algebroids and Lie pseudoalgebras,
{\em Bull. London Math. Soc. } {\bf 27} (1995), 97--147.  
\bibitem{ma-xu:lie}
Mackenzie, K.C.H. and Xu, P.,
Lie bialgebroids and Poisson groupoids,
{\em Duke Math. J.} {\bf 73} (1994). 415-452.
\bibitem{ro:courant}
Roytenberg, D., Courant algebroids, derived brackets and
even symplectic supermanifolds, Ph.D. thesis,
University of California, Berkeley (1999), preprint math.DG/9910078.
\bibitem{ro-we:courant} 
Roytenberg, D., and Weinstein, A., Courant algebroids
and strongly homotopy Lie algebras, {\em Lett. Math. Phys.} {\bf 46}
(1998), 81--93.
\bibitem{va:lectures} 
Vaisman, I., {\em Lectures on the Geometry of Poisson Manifolds},
Birkh\"{a}user, Basel, 1994.
\end{thebibliography}
\end{document}